\documentclass{tac}
\usepackage{xy}
\usepackage{MnSymbol}
\usepackage{tikz}

\xyoption{all}

\def\mathcaldef#1{\expandafter\def\csname#1\endcsname{{\cal#1}}}

\def\"{``}
\def\q{\quad}
\def\qq{\quad\quad}

\def\imp{\Rightarrow}

\def\op{^{\rm op}}
\def\ov{\overline}

\def\tm{\times}

\def\si{\sigma}
\def\t{_\blacktriangleright}

\def\cd{\ldots}

\mathcaldef{C}
\mathcaldef{D}
\mathcaldef{M}
\mathcaldef{N}
\mathcaldef{L}
\mathcaldef{B}
\mathcaldef{E}
\mathcaldef{I}
\mathcaldef{G}

\mathbfdef{Set}
\mathbfdef{Set_f}
\mathbfdef{cMon}
\mathbfdef{Mon}
\mathbfdef{Mlt}
\mathbfdef{fpMlt}
\mathbfdef{sMlt}
\mathbfdef{Rep}
\mathbfdef{fpRep}
\mathbfdef{sRep}
\mathbfdef{Seq}
\mathbfdef{fpSeq}
\mathbfdef{Cat}
\mathbfdef{fpCat}
\mathbfdef{Sum}
\mathbfdef{preAdd}
\mathbfdef{Add}
\mathbfdef{n}
\mathbfdef{m}

\mathrmdef{id}

\def\obj{{\rm obj\,}}
\def\sq{{\rm sq}}
\def\sp{{\rm sp}}
\def\pb{{\rm pb}}
\def\tr{{\rm tr}}
\def\sptr{{\rm sptr}}

\newtheorem{prop}{Proposition}

\let\pf\proof
\let\epf\endproof
\def\eq{\begin{equation}}
\def\eeq{\end{equation}}

\author{Claudio Pisani}
\address{via Saluzzo 67,\\ 10125 Torino, Italy.}
\title{Fibered multicategory theory}
\copyrightyear{2022}
\eaddress{pisclau@yahoo.it}

\begin{document}

\maketitle

\begin{abstract}

Given a fibration in groupoids $d:\D\to\I$, we define a fibered multicategory as a particular functor $p:\M\to\I$,
where $\M$ has the same objects as $\D$, and its arrows $a:X\to Y$ should be thought of as 
families of arrows in the multicategory, indexed by $pY$.
The key axiom extends the reindexing of objects, given by $d$, to a reindexing of arrows in $\M$
along pullback squares in $\I$.
When $\D$ is included in $\M$, in an appropriate sense, one gets again fibered categories. 
In this context, cartesian fibered multicategories are defined and studied in a natural way.

\end{abstract}


\section{Introduction}
\label{intro}

In the last decades, generalized multicategories have been the object of renewed interest 
(see for instance \cite{leinster}, \cite{hermida2} and \cite{shulman}).
The focus is usually set on the monad that gives the generalized domains of arrows, and on 
the appropriate environment for the monad itself.
Our emphasis here lies rather on the fibered aspects of a symmetric multicategory $\M$, which derive from 
considering {\em families} of arrows in $\M$ as the arrows of a {\em category} over finite sets.
The resulting abstract notion of \"fibered multicategory" $p:\M\to\I\,$ generalizes that of fibered category
in that cartesian arrows (giving \"reindexing of objects") are not in principle included in $\M$.

\begin{remark}
Although there is a potential conflict of terminology with the concept of
\"fibration of multicategories" studied in \cite{hermida2} and in related works, 
note that notion refers instead to a particular kind of {\em morphism of multicategories}. 
So, the context should mitigate the risk of confusion. 
\end{remark} 

The present theory has developed from the attempt to come to terms with two aspects 
of symmetric multicategories.
On the one hand, in the classical definition one considers finite {\em sequences} of objects 
as the domain of arrows. However, this seems a somewhat unnatural restriction 
(analogous, and in fact related, to the choice of working in skeletal categories)
and one may want to allow arbitrary finite families of objects as domain.
This step was in fact already considered by \cite{beilinson} (\"pseudo-tensor categories") 
and by \cite{leinster} (\"fat multicategories"), though those definitions seem somewhat {\em ad hoc},
especially in the way they treat identities.

On the other hand, there is the natural idea of associating, to a symmetric multicategory $\M$, 
the ordinary category of finite families (or sequences) of arrows in $\M$. 
In the case of operads, this construction dates back at least to \cite{may} where,
for some reasons, these \"categories of operators" are then considered as categories 
over pointed finite sets (rather than on finite sets), somehow rendering less sharp the concept.
The construction is also related to the free strict monoidal category on $\M$
(see for instance \cite{leinster}).

Embracing wholeheartedly both these perspectives, we are led to consider the category $\M'$ of finite families 
of objects and of arrows in $\M$, with its projection functor $p:\M'\to\Set_f$.
Thus, an arrow $a:X\to Y$ in $\M'$ (of \"shape" $pa: I \to J$ in $\Set_f$)
is a family of arrows in $\M$, indexed by $J = p Y$,  
and the \"multicategorical" composition in $\M$ is encoded in the \"categorical" composition in $\M'$.
The other piece of structure on $\M'$ to take in account is the usual reindexing for its objects $X = x_i \,\,\, (i\in I)$,
which can be encoded by a discrete fibration $d:\D \to \Set_f$, where $\D$ has the same objects as $\M'$.

Now, in order to develop a sensible theory of multicategories relative to a general base category $\I$,
it is natural to look for the relevant relations between the arrows of $\M'$ and those of $\D$.
Among them, it seems central the possibility of reindexing arrows in $\M'$ along pullbacks in $\Set_f$: 
{
given $a:X\to Y$ in $\M'$, given a pullback $(pa)f' = fg$ in $\Set_f$ whose right side is $pa$, and given
$d$-cartesian liftings of its bottom and top sides $f$ and $f'$, there is a unique \"special" lifting $f^*a$ 
of the left side.
}
\[
\xymatrix@R=4pc@C=4pc{
f^*X \ar@{..>}[d]_{f^*a}\ar@{-->}[r]^{f'}  & X \ar[d]^a \\
f^*Y \ar@{-->}[r]_f                    & Y      }
\xymatrix@R=2pc@C=2pc{
 \\ & \ar@{~>}[r]^p &     } \qq
\xymatrix@R=1.3pc@C=1.3pc{
K \ar[dd]_g\ar[rr]^{f'}  & & pX \ar[dd]^{pa} \\
& {\rm pb} & \\
L  \ar[rr]_f                 & & pY       }
\]

\begin{remark}
It should be noted that the notation $f^*a$ is somewhat misleading, since that arrow 
depends in fact on the whole pullback square. 
Note also that throughout the paper we use dashed arrows for the maps in $\D$ 
(and often will name them as the map they lift), dotted arrows for the uniquely determined ones, 
while the wavy arrow refers to the projection to the base $\I$.
\end{remark}

For instance, if $1$ is a terminal set and the bottom side of the pullback, $j:1\to J$, picks an element of the 
indexing set $J= pY$, that \"special" lifting is a corresponding single arrow of the family $a:X\to Y$.
A general mapping $f:L\to J$ gives a new family $f^*a: f^*X\to f^*Y$, containing repeated copies 
of some the arrows in the family $a$ (if $f$ is not injective) and deleting others (if it is not surjective).
On the other hand, these \"copies" have in general a different indexing.
Indeed, if $f:L\to J$ is a bijective mapping, $f^*a: f^*X\to f^*Y$ is the \"same" as $a:X\to Y$ apart from
a bijective reindexing of the objects involved (both in $X$ and in $Y$). 
Nonetheless, the indexing is necessary in order to have a non-ambiguous composition of families. 
Thus, for instance, a bijective reindexing of an identity in $\M'$ is not in general an identity (if $f'\neq f$).

Now, it is not difficult to envisage the corresponding abstract setting.
First, we have a category $\I$ with pullbacks and two categories over it, $d:\D\to\I$ and $p:\M\to\I$,
sharing the same objects.
These data allow us to consider the double category of squares with \"horizontal" arrows in $\D$
and \"vertical" arrows in $\M$:
\[
\xymatrix@R=4pc@C=4pc{
U \ar[d]_b\ar@{-->}[r]^{f'}  & X \ar[d]^a \\
V \ar@{-->}[r]_f                    & Y      }
\qq\qq\qq\qq
\xymatrix@R=1.3pc@C=1.3pc{
U \ar[dd]_b\ar@{-->}[rr]^{f'}  & & X \ar[dd]^a \\ 
& \# & \\
V \ar@{-->}[rr]_f                 & & Y      }
\]
Next we assume that there is a sub double category of \"special" squares (marked, as above,
by the symbol $\#$), and that these are projected to pullback squares in $\I$.
Lastly, we assume the following fibration-like properties:
\begin{enumerate}
\item 
$d:\D\to\I$ is a groupoid fibration. (We do not assume it discrete for not ruling out significant instances,
such as the \"self indexing" of the multicategory $\Set$, see section \ref{setset}).
\item
Special squares form a discrete fibration over pullbacks with chosen $d$-liftings of the top
and the bottom sides (see Section \ref{fm} for the precise definition):
\end{enumerate}

\eq   \label{int3}
\xymatrix@R=4pc@C=4pc{
U \ar@{-->}[r]^{f'}  & X \ar[d]^a \\
V \ar@{-->}[r]_f                    & Y      }
\qq\qq\qq\qq
\xymatrix@R=1.3pc@C=1.3pc{
U \ar@{..>}[dd]_b\ar@{-->}[rr]^{f'}  & & X \ar[dd]^a \\ 
& \# & \\
V \ar@{-->}[rr]_f                 & & Y      }
\eeq
\[
\xymatrix@R=1.3pc@C=1.3pc{
K \ar[dd]_g\ar[rr]^{f'}  & & pX \ar[dd]^{pa} \\
& {\rm pb} & \\
L  \ar[rr]_f                 & & pY       }
\]

Classical symmetric multicategories (albeit in their \"fat" form) arise in this conceptual frame when 
\begin{enumerate}
\item 
$\I=\Set_f$ and $d:\D\to\Set_f$ is the discrete fibration represented 
(via the inclusion $\Set_f\to\Set$) by a set $\M_0$;
\item
the following \"extensivity" condition holds: 
for any finite family of pullbacks in $\Set_f$ along injections 
of a finite sum $i_k:I_k\to I$, 
\[
\xymatrix@R=1.3pc@C=1.3pc{
J_1 \ar[rr]\ar[dd]_{f_1}  && J \ar[dd]^f && J_2 \ar[dd]^{f_2}\ar[ll]  \\
& {\rm pb} && {\rm pb} \\
I_1 \ar[rr]_{i_1}       && I   && I_2 \ar[ll]^{i_2}       }
\]
and for any lifting $a_k:X_k\to Y_k$ of its external vertical sides, and $d$-liftings of the horizontal sides,
there is a unique lifting $a:X\to Y$ of $f$ forming special squares:
\eq   \label{int5}
\xymatrix@R=1.3pc@C=1.3pc{
X_1 \ar@{-->}[rr]\ar[dd]_{a_1}  && X \ar@{..>}[dd] && X_2 \ar[dd]^{a_2}\ar@{-->}[ll] \\
& \# && \# \\
Y_1 \ar@{-->}[rr]_{i_1}       && Y   && Y_2 \ar@{-->}[ll]^{i_2}       }
\eeq
\end{enumerate}
This condition assures that \"abstract families" $a:X\to Y$ in $\M$ \"are" in fact families 
of \"single arrows", that is of arrows $s:X\to Y$ with $pY$ terminal.
It is then easy to see that both the action on arrows of permutations of the objects in the domain
(which are a particular case of bijective reindexing) and their relationships with composition 
are direct consequences of the reindexing axiom (\ref{int3}) of arrows along pullbacks. 

\subsection{Outline of the paper}\qq

In Section \ref{fm}, fibered multicategories are defined and confronted with 
classical symmetric multicategories.
We also present the related concepts of fibered monoid and of monoid in a fibered multicategory.
In sections \ref{rep} and \ref{unipro}, we define representability and universal products.
In section \ref{comm}, we show how the concept of commuting endomorphisms
also finds a natural place in this abstract framework.

In Section \ref{fc} we consider the connections of fibered multicategories with ordinary fibered categories.
Of course, the latter are included in the former: given a fibered category $p:\C\to\I$, we get a
fibered multicategory by considering it paired with the restriction $d:\D\to\I$ of $p$ to the $p$-cartesian arrows;
the special squares are then simply the commutative squares in $\C$ (with cartesian horizontal arrows) 
over pullbacks in $\I$.
If $\C$ is an ordinary category, and $p:\C\to\Set_f$ is the usual family fibration,
the correponding fibered multicategory is associated to the \"sequential" multicategory $\C\t$ of families
of concurrent arrows (see \cite{pisani}).

Conversely, we will see that when $\D$ can be included in $\M$ in a precise sense, 
then $p:\M\to\I\,\,$ is in fact a fibered category.
In this context, the Grothendieck construction can be seen as the free way to transform 
a fibered unary multicategory $p:\M\to\I$ 
(that is, one in which all the arrows of $\M$ are projected to isomorphisms) into a fibered category.
The fact that the construction yields an equivalence reflects (in the more general context) 
the elementary fact that unary multicategories are equivalent to categories (or to sequential multicategories).

In Section \ref{cm}, we show how in this abstract framework there is an effective notion of cartesian structure.
A fibered cartesian multicategory has, in addition to the special squares, also \"special triangles"
giving \"covariant reindexing" of arrows.
These special triangles form a discrete opfibration over commutative triangles in $\I$
(with a chosen $d$-lifting of its top side)
and are related to special squares by some sort of Beck-Chevalley and Frobenius laws.
As expected, in a cartesian multicategory representability (by a tensor product) becomes an
\"algebraic" or \"absolute" property and coincides with the existence of \"universal" or \"cartesian" products.


\section{Fibered multicategories}
\label{fm}

In order to define smoothly fibered multicategories, it seems appropriate to consider explicitly
some double categories that arise in a natural way, as in the following constructions.

\subsection{Construction A}

Let $d:\D\to\I$ be any functor, and let $\pb(d,d)$ the subcategory of the comma category $(d,d)$
whose arrows are pullback squares in $\I$.
Thus, an object of $\pb(d,d)$ consists of a pair of objects $X,Y\in\D$ and an arrow $l:dX\to dY$ in $\I$;
an arrow from $l:dX\to dY$ to $h:dX'\to dY'$ consists of arrows $f:X\to X'$ and $g:Y\to Y'$ in $\D$,
such that their images form a pullback square $h(df) = (dg)l$ in $\I$:

\[
\xymatrix@R=4pc@C=4pc{
X \ar[r]^f  & X' \\
Y \ar[r]^g                   & Y'      }
\xymatrix@R=2pc@C=2pc{
 \\ & \ar@{~>}[r]^d &      }    \qq
\xymatrix@R=1.3pc@C=1.3pc{
dX \ar[dd]_l\ar[rr]^{df}  & & dX' \ar[dd]^h \\
& \pb & \\
dY \ar[rr]_{dg}                & & dY'       }
\]

Since the arrows of $\pb(d,d)$ are squares which can also be composed vertically in the obvious way,
$\pb(d,d)$ is in fact a double category.

\begin{remark}
Note that if $d:\D\to\I$ is a fibration in groupoids (that is a fibration in which all arrows are cartesian, 
so that the fibers are groupoids) then also the obvious functor $\pb(d,d)\to\pb\,\I$ is such.
\end{remark}

\subsection{Construction B}

Let $\D$ and $\M$ be categories with the same objects. 
We can consider the double category $\sq(\M,\D)$ which has $\M$ as horizontal category,
$\D$ as vertical category, and all squares as arrows.
If furthermore there are functors to a third category $d:\D\to\I$ and $p:\M\to\I$,
we can consider the sub double category $\pb(\M,\D)$ of those squares in $\sq(\M,\D)$
projected to pullback squares in $\I$.

\[
\xymatrix@R=4pc@C=4pc{
X \ar[r]^f\ar[d]_{a}        & X'\ar[d]^{b}  \\
Y \ar[r]_g      & Y'      }
\xymatrix@R=2pc@C=2pc{
 \\ & \ar@{~>}[r] &      } \qq
\xymatrix@R=1.3pc@C=1.3pc{
dX \ar[dd]_{pa}\ar[rr]^{df}  & & dX' \ar[dd]^{pb} \\
& {\rm pb} & \\
dY \ar[rr]_{dg}                & & dY'       }
\]

\begin{remark}
There is an obvious double functor $\pb(\M,\D)\to\pb(d,d)$ over $\pb\,\I$.
\end{remark}

\subsection{Fibered multicategories}

We are now in a position to define the main notion of the present paper.

\begin{definition}
A {\em fibered multicategory} consists of
\begin{enumerate}
\item
A category $\I$ with pullbacks.
\item
A fibration in groupoids $d:\D\to\I$.
\item
A category $\M$ such that $\obj\M = \obj\D$, and a functor $p:\M\to\I$ which coincides with $d$ on objects.
\item \label{key}
A sub double category $\sp(\M,\D)$ of $\pb(\M,\D)$, such that the obvious functor 
\[ \sp(\M,\D) \to \pb(d,d) \] 
is a discrete fibration, with respect to horizontal arrows.
\end{enumerate}

A {\em morphism} $F:\M \to \M'$ of fibered multicategories over $\I$, consists of a 
pair of functors over $\I$, $F_\D:\D\to\D'$ and $F_\M:\M\to\M'$ which coincide on objects 
and which take special squares to special squares.
\end{definition}

We often indicate the fibered muticategory $(d, p, \sp(\M,\D))$ simply by $\M$, 
leaving the other data implicit.
Squares in $\sp(\M,\D)$ are called \"special" squares and are graphically marked with the symbol $\#$.

As sketched in the introduction, the key axiom (namely, point \ref{key} in the definition) 
says that the following \"unique special lifting" property holds:
{\em given an arrow $a:X\to Y$ in $\M$, a pullback $(pa)f' = fg$ in $\I$, 
and a pair of $d$-liftings of $f'$ and $f$ (with codomain $X$ and $Y$),
there is a unique special square over the given pullback extending these data}:

\eq   \label{fm3}
\xymatrix@R=4pc@C=4pc{
U \ar@{-->}[r]^{f'}  & X \ar[d]^a \\
V \ar@{-->}[r]_f                    & Y      }
\qq\qq\qq\qq
\xymatrix@R=1.3pc@C=1.3pc{
U \ar@{..>}[dd]_b\ar@{-->}[rr]^{f'}  & & X \ar[dd]^a \\ 
& \# & \\
V \ar@{-->}[rr]_f                 & & Y      }
\eeq
\[
\xymatrix@R=1.3pc@C=1.3pc{
K \ar[dd]_g\ar[rr]^{f'}  & & pX \ar[dd]^{pa} \\
& {\rm pb} & \\
L  \ar[rr]_f                & & pY       }
\]

\subsection{Acting on maps}
\label{action}

There are various ways in which one can act on maps in $\M$ exploiting the unique special lifting property.
For instance, since $\sp(\M,\D)\to\pb(d,d)$ is a discrete fibration, and both $\pb(d,d)\to\pb\,\I$ and $\pb\,\I\to\I$
are fibered in groupoids, their composition $\sp(\M,\D)\to\I$ (taking arrows
in $\M$ and special squares to their codomain) is itself fibered in groupoids.
Thus, an arrow $a:X\to Y$ in $\M$ can be pulled back (uniquely up to a unique special square) 
along an arrow $f:I\to pY$ in $\I$:
\[
\xymatrix@R=1.3pc@C=1.3pc{
U \ar[dd]_{f^*a}  & & X \ar[dd]^a \\
&  & \\
V                   & & Y   \\
I  \ar[rr]^f               & & pY  }
\]

On the other hand, if we restrict ourselves to pullbacks with an identity bottom side,
then $X\mapsto \M(X,Y)$ extends, for any $Y\in\M$, to a functor $\si_Y:\G\op\to\Set$, 
where $\G$ is the groupoid of isomorphisms in $\D$.
\eq   \label{fm4}
\xymatrix@R=4pc@C=4pc{
U \ar@{-->}[r]^f  & X \ar[d]^a \\
Y \ar@{-->}[r]_\id                   & Y      }
\qq\qq\qq\qq
\xymatrix@R=1.3pc@C=1.3pc{
U \ar@{..>}[dd]_{f^*a}\ar@{-->}[rr]^f  & & X \ar[dd]^a \\ 
& \# & \\
Y \ar@{-->}[rr]_\id                 & & Y      }
\eeq
\[
\xymatrix@R=1.3pc@C=1.3pc{
K \ar[dd]\ar[rr]^f  & & pX \ar[dd]^{pa} \\
& {\rm pb} & \\
pY  \ar[rr]_\id                 & & pY       }
\]
The functoriality of the $\si_Y$ follows from the assumption that special squares compose horizontally:
\eq   
\xymatrix@R=1.3pc@C=1.3pc{
U \ar@{-->}[rr]\ar[dd]  && V \ar[dd]\ar@{-->}[rr] && X \ar[dd] \\
& \# && \# \\
Y \ar@{-->}[rr]_\id      && Y\ar@{-->}[rr]_\id   &&   Y     }
\eeq

If $d$ is itself a {\em discrete} fibration, then, for each $X$ and $Y$, the automorphism
group on $p X$ in $\I$ acts on the set $\M(X,Y)$, like in the classical definition
of symmetric multicategories.
The fact that special squares compose vertically
\eq   
\xymatrix@R=1.3pc@C=1.3pc{
V \ar@{-->}[rr]\ar[dd]   && W \ar[dd]  \\
& \# \\
U \ar@{-->}[rr]\ar[dd]   && X \ar[dd]  \\
& \# \\
Y \ar@{-->}[rr]_\id        && Y   }
\eeq
encodes the classical \"covariance" conditions of the action with respect to composition 
in a symmetric multicategory.

\subsection{Standard multicategories}
\label{standard}

Now we sketch how our abstraction subsumes the classical notion of symmetric multicategory, 
or better, the variant with general indexing of domains of arrows, as argued in the introduction.
In fact, if $\I = \Set_f$ and $d$ is the discrete fibration \"represented" by a set $\M_0$
then, assuming also the extensivity axiom (\ref{int5}), we get our version of the  
\"fat multicategories" of \cite{leinster} or of the \"pseudo-tensor categories" of \cite{beilinson}. 

\begin{definition}
Given a set $\M_0\in\Set$, denote by $\D$ the corresponding comma category  $i/\M_0$
(where $i$ is the inclusion $\Set_f\to\Set$) with its projection $d : \D \to \Set_f$.
A {\em standard multicategory}  with \"set of objects" $\M_0$ is a fibered multicategory 
with $\D$ as above and satisfying condition (\ref{int5}).

A {\em standard morphism} $F:\M \to \M'$ between standard multicategories is a morphism 
such that the component $F_\D:\D\to\D'$ is induced by a mapping $F_0:\M_0\to\M'_0$.
\end{definition}

Thus, in a standard multicategory, the objects are finite families of elements in $\M_0$, 
and the reindexing of objects is the usual one.

\begin{remark}
If $I\in\Set_f$ has just one element, for any $x\in\M_0$ there is an identity $\id_x:x\to x$ 
over the identity of $I$.
If $J$ is {\em another} terminal set, one can reindex the domain of $\id_x$ to obtain an arrow in $\M$ 
which {\em is not} an identity, since domain and codomain are different indexing of $x$, 
and so are different objects in $\M$:
\[
\xymatrix@R=1.3pc@C=1.3pc{
x_J \ar@{..>}[dd]\ar@{-->}[rr]  & & x_I \ar[dd]^\id \\ 
& \# & \\
x_I \ar@{-->}[rr]^\id                 & & x_I      }
\xymatrix@R=2pc@C=2pc{
 \\ & \ar@{~>}[r] &      } \qq
\xymatrix@R=1.3pc@C=1.3pc{
J \ar[dd]\ar[rr]  & & I \ar[dd] \\
& {\rm pb} & \\
I \ar[rr]                & & I       }
\]
Such \"pseudo-identities" are called \"identity maps" in the definition of fat multicategories in \cite{leinster},
and are exploited therein to get the reindexing of domain of arrows.
\end{remark}

Now, exploiting axiom (\ref{int5}) which gives the equivalence between arrows in $\M$ 
and families of \"single" arrows $a:X\to Y$ with $pY$ terminal, 
it is easy to see that our standard multicategories are essentially the same as fat multicategories.
Thus, following \cite{leinster}, one also shows that the category of standard multicategories 
is equivalent to the category of classical symmetric multicategories; 
essentially, the latter are the skeletal version of the former,
with respect to the skeleton ${\bf N} \subseteq \Set_f$  formed by the sets 
$\n =\{1, 2 \cdots, n\}\,\,$ (including ${\bf 0} = \emptyset$).

\subsection{The multicategory of sets indexed over sets}
\label{setset}

Now we see an example of fibered multicategory in which the reindexing fibration is not discrete.
Let $\I = \Set$ and let $d:\D\to\Set$ be the usual codomain fibration restricted to pullback squares.
Given a mapping $f:I\to J$ in $\Set$ and objects $q:X\to I$ and $r:Y\to J$ in $\D$,
a map in $\M$ over $f$ consists of a mapping which takes a family of elements in $X$ 
\[ x_{i j} \q (j\in J) \q (i\in f^{-1}j)   \]
such that $q(x_{i j}) = i$, to a family $\q y_j \q (j\in J)$ of elements in $Y$  such that $r(y_j) = j$.
Now, suppose that in the diagram below the wavy arrow on the front represents a map over $f$ as
just described, and suppose that the bottom, left and right faces are pullbacks.
We leave it to the reader to envisage the unique sensible way to obtain the wavy arrow on the back, 
in order to define special squares.
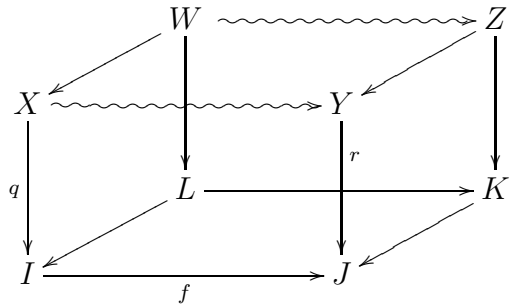
\[
\xymatrix@R=1.5pc@C=1.5pc{
&&& W \ar[dd]\ar[lld]  \ar@{~>}[rrrr] &&&& Z\ar[dd]\ar[lld]  \\
&X \ar[dd]_q \ar@{~>}[rrrr] &&&& Y\ar[dd]^(.3)r \\ 
&& & L \ar[lld]\ar[rrrr] &&&& K \ar[lld] \\
& I \ar[rrrr]_f  &&&& J  &      }   
\]
Of course, this is nothing but the self indexed version of the classical multicategory $\Set$,
whose arrows are many variable mappings.

\subsection{Representable and stably representable multicategories}
\label{rep}

A fibered multicategory $\M = (d, p, \sp(\M,\D))$ is {\em representable} if $p$ is an opfibration. 
Of course, in the standard case we get the usual notion (see \cite{hermida}) which
gives a universal characterization of tensor products.
We say that an arrow is {\em stably opcartesian} if all its reindexing are opcartesian, and
that $\M$ is {\em stably representable} if it is representable and the opcartesian arrows are 
stably opcartesian.

\subsection{Universal products}
\label{unipro}

If $X\in\M$ and $f:pX\to J$ is a map in $\I$, a {\em universal product} for $X$ along $f$ 
is an object $P\in\M$ over $J$ along with a vertical map $\pi:f^*P\to X$ 
\eq  \label{unipr}
\xymatrix@R=1.3pc@C=1.3pc{
f^*P \ar[dd]_\pi\ar@{-->}[ddrr]^f  & &  \\ 
& &  \\
X                 & & P      } 
\xymatrix@R=2pc@C=2pc{
 \\ & \ar@{~>}[r] &      } \qq
\xymatrix@R=1.3pc@C=1.3pc{
pX \ar[dd]_\id\ar[ddrr]^f  & &  \\ 
& &  \\
pX                 & & J      } 
\eeq
which has the following universal property:
for any pullback $hf' = fh'$ in $\I$, any $d$-lifting or its top side $f'^*Q\to Q$,
and any arrow $\rho: f'^*Q\to X$ over $h'$, there exists a unique completion 
to a special square $(t,t',f,f')$ such that $\rho = \pi t'$.
\[
\xymatrix@R=1.3pc@C=1.3pc{
& f'^*Q \ar@{..>}[ddl]_{t'} \ar@/^/[ddddl]^(0.3)\rho \ar@{-->}[ddrr]^{f'}&& \\ \\
f^*P \ar[dd]_\pi\ar@{-->}[ddrr]^f  & & \# & Q \ar@{..>}[ddl]^t \\ 
& &  \\
X                 & & P      } 
\xymatrix@R=2pc@C=2pc{
 \\ \\ & \ar@{~>}[r] &      } \qq
\xymatrix@R=1.3pc@C=1.3pc{
& L \ar[ddl]_{h'} \ar@/^/[ddddl]^(0.3){h'} \ar[ddrr]^{f'}&& \\ \\
pX \ar[dd]_\id\ar[ddrr]^f  & & \pb & K \ar[ddl]^h \\ 
& &  \\
pX                 & & J      } 
\]
In the standard case, it is easy to see that a universal product for $X = x_i \,\, (i\in I)$
along $f:I\to 1$, is an object $x\in\M_0$ together with projections $\pi_i:x\to x_i$
such that for any $Y\in\M$ and any family of arrows $\rho_i:Y\to x_i$, there is a unique $t:Y\to x$
with $\pi_i t = \rho_i \,\, (i\in I)$.

If the in above definition of universal product we restrict to pullbacks in $\I$ of the form $\id f = f \id$,
so that $\rho$ is vertical, we get in the standard case the weaker condition in which $pY$ is also terminal
(that is universality holds only for unary arrows $\rho_i:y\to x_i$).

\subsection{Fibered monoids}
\label{fibmon}

\qq A {\em fibered monoid} is a representable discrete multicategory $\M = (d, p, \sp(\M,\D))$,
that is, one such that $p$ is a {\em discrete} opfibration.

If $d$ is also discrete, the notion simplifies further, since the special lifting is forced by the underlying pullback.
In this case, a monoid over $\I$ consists of a discrete fibration $d$ and a discrete opfibration $p$ 
with the same objects, such that for any pullback in $\I$, and any lifting $X$ of its top-right object, 
the two possible paths from $X$ over the given pullback, through vertical $p$-liftings and horizontal $d$-liftings, 
end up with the same lifting of the bottom-left object. 

Of course, in the standard case, one obtains the \"unbiased" and \"fat" notion of commutative monoid
on a set $\M_0$.
In section \ref{monin} below, we will define monoids {\em in } a fibered multicategory $\M$,
which instead need not to be commutative.

\subsection{Commutative arrows and commuting endomorphisms}
\label{comm}

Two concepts that can be formulated in symmetric multicategories are those of {\em commutative arrows}
and of {\em commuting endomorphisms}.
In our context, we have seen in section \ref{action} that for any $Y\in\M$ we have a functor 
$\si_Y:\G\op\to\Set$, $X\mapsto \M(X,Y)$, where $\G$ is the groupoid of isomorphisms in $\D$. 
In particular, any arrow in $\M$ has its automorphism group, giving its symmetries. 
These can be seen as a measure of its degree of commutativity.

As for the second notion, we begin by defining endomorphisms.
\begin{definition}
\label{endo}
An {\em endomorphism} is an arrow $a:X\to Y$ in $\M$ along with a map $t$ in $\D$ parallel to it
with $d t = p a$.
\end{definition}
\[
\xymatrix@R=2pc@C=2pc{
X\ar@<-0.5ex>[rr]_a \ar@<0.5ex>@{-->}[rr]^t &&Y     }
\]
Note that in general $t$ is not uniquely determined by $a$, while of course this is the case if $d$ is discrete.
It is easy to see that, in the standard case, an endomorphism $a:X\to Y$ is a family of 
endoarrows $a_i : x_i,\cd ,x_i \to x_i$.

Now, suppose that we have two endomorphisms $a,t : X\to Y$ and $a',t' : X'\to Y$ with the same codomain,
a pullback in $\I$ completing their images and a $d$-lifted completion of that pullback:
\[
\xymatrix@R=1.6pc@C=1.5pc{
 & & X \ar@<-0.5ex>[dd]_a \ar@<0.5ex>@{-->}[dd]^t \\ 
 \\
X' \ar@<-0.5ex>@{-->}[rr]_{t'}  \ar@<0.5ex>[rr]^{a'}      & & Y      }
\qq\qq
\xymatrix@R=1.3pc@C=1.3pc{
K \ar[dd]\ar[rr]  & & I \ar[dd] \\
& {\rm pb} & \\
I' \ar[rr]                & & J       }
\qq\qq
\xymatrix@R=1.6pc@C=1.5pc{
Z \ar@{-->}[dd]\ar@{-->}[rr]  & & X \ar@{-->}[dd]^t \\
&& \\
X' \ar@{-->}[rr]_{t'}                 & & Y       }
\]
By the reindexing axiom (applied once \"horizontally" and once \"vertically") 
we get then the two dotted maps below:
\eq
\label{com}
\xymatrix@R=1.6pc@C=1.5pc{
Z \ar@<-0.5ex>@{..>}[dd]_b\ar@<0.5ex>@{-->}[dd] \ar@<-0.5ex>@{-->}[rr]\ar@<0.5ex>@{..>}[rr]^{b'} & & 
X \ar@<-0.5ex>[dd]_a \ar@<0.5ex>@{-->}[dd]^t \\ 
& \# &  \\
X' \ar@<-0.5ex>@{-->}[rr]_{t'}  \ar@<0.5ex>[rr]^{a'}      & & Y      }
\eeq
We say that the endomorphisms $a$ and $a'$ {\em commute} if the above diagram 
commutes in $\M$: $ab' =a'b$. The notion does not depend on the choice of the pullback
and of its $d$-lifting. Indeed, suppose that $ab' =a'b$ and consider another such a choice; 
one has then a mediating $d$-arrow $u$, as below:
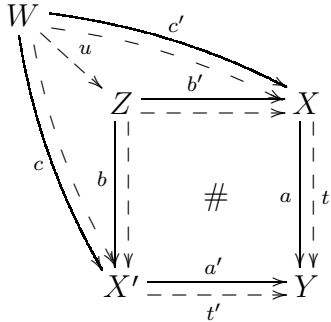
\[
\xymatrix@R=1.6pc@C=1.5pc{
W \ar@{-->}[dr]^u \ar@<-0.5ex>@/_/[dddr]_c\ar@<0.5ex>@/_/@{-->}[dddr] 
\ar@<0.5ex>@/^/[drrr]^{c'}\ar@<-0.5ex>@/^/@{-->}[drrr] \\
& Z \ar@<-0.5ex>[dd]_(0.4)b\ar@<0.5ex>@{-->}[dd] \ar@<-0.5ex>@{-->}[rr]\ar@<0.5ex>[rr]^(0.4){b'} & & 
X \ar@<-0.5ex>[dd]_a \ar@<0.5ex>@{-->}[dd]^t \\ 
&& \# &  \\
& X' \ar@<-0.5ex>@{-->}[rr]_{t'}  \ar@<0.5ex>[rr]^{a'}      & & Y      }
\]
By the fibrational properties of special squares, the triangles $(c,b,u)$ and $(c',b',u)$ 
are also special squares (with an identity side). 
By composition, also $(a'c,a'b,u)$ and $(ac',ab',u)$ are special.
Since $ab' =a'b$ by hypothesis, $ac'$ and $a'c$ are the reindexing of the same map 
along the same square in $\pb(d,d)$, so that they are the same.

Of course, again, in the standard case one gets the usual notion.
Note indeed that for single arrows $a:X\to Y$ and $a':X'\to Y$ (with $pY$ terminal), in the square (\ref{com})
above, $Z$ is over a product $pX\times pX'$, and the two paths along the square indicate
the two possible ways of combining $a$ and $a'$ to get a single arrow.

\subsection{Monoids in a fibered multicategory}
\label{monin}

A {\em monoid in} $\M = (d, p, \sp(\M,\D))$ is a morphism $m:1 \to \M$ from the terminal multicategory 
over $\I$.
So it amounts to a section of $d$ and a section of $p$ that coincide on objects and such that 
pullbacks are lifted to special squares. 
We may say that a section of $d$ is an {\em object in } $\M$, so that any monoid 
in $\M$ has an underlying object. Indeed, in the standard case the sections of $d$ correspond
to the elements $x\in\M_0$ (since $\M_0$ is the limit of the functor $\D\to\Set_f\to\Set$).

Thus, a monoid $m$ in $\M$ gives, for any map $f:I\to J$ in $\I$, an endomorphism $(t_f,a_f)$ over it
in such a way that these endomorphisms are closed with respect to composition and identities.
\[
\xymatrix@R=2pc@C=2pc{
X\ar@<-0.5ex>[rr]_{a_f} \ar@<0.5ex>@{-->}[rr]^{t_f} &&Y     }
\]

A {\em morphism of monoids} $\alpha : m \to m'$ consists of vertical maps in $\M$, 
\[
\alpha_I:X\to X' \qq (I\in\I)
\] 
such that, for any $f:I\to J$, the square 
\[
\xymatrix@R=1.6pc@C=1.5pc{
X \ar[dd]_{\alpha_I} \ar@<-0.5ex>@{-->}[rr]_{t_f}\ar@<0.5ex>[rr]^{a_f} & & 
Y \ar[dd]^{\alpha_J} \\ 
& \# &  \\
X' \ar@<-0.5ex>@{-->}[rr]_{t'_f}  \ar@<0.5ex>[rr]^{a'_f}      & & Y'      }
\]
is special (with respect to $t$ and $t'$) and commutative (with respect to $a$ and $a'$).
Thus, we may say that the $\alpha_I$ give a \"special" natural transformation between the 
$p$-sections $m$ and $m'$.

\subsection{Eckmann-Hilton arguments}
\label{hilton}

In such an unbiased definition of monoid, there is no preferred \"associative" or \"identity" law.
For instance, for any $I\to J\to 1$ in $\I$, the corresponding commutative triangle of endomorphisms
can be seen as an associative law, while for any section-retraction pair in $\I$, such as $1\to I\to 1$,
the corresponding triangle can be seen as an identity law.

Relevant to the Eckmann-Hilton arguments in the present context seems to be the \"diagonal identity law",
associated to the diagonal-projection pair for an object $I\in\I$.
\[
\xymatrix@R=2pc@C=3pc{
I\ar[r]_\Delta \ar@/^1.5pc/[rr]^\id & I\times I \ar[r]_\pi  &  I  }
\]

Accordingly, we say that two monoids $m$ and $m'$ in $\M$ have the same identity over $I\in\I$
if $\Delta: I \to I\times I$ is lifted to the same endomorphism $a_\Delta$ by $m$ and $m'$.
\begin{proposition}
Suppose that the monoids $m$ and $m'$ have the same identity over $I$,
and that the liftings $a_!$ and $a'_!$ of $\,\,! : I \to 1$ commute.
Then $a_! = a'_!$.
\end{proposition}
\pf
Indeed, in the left hand diagram below, which is over the right one,
the square  and the triangles commute by hypothesis.
\[
\xymatrix@R=1.6pc@C=1.5pc{
X \ar[dr]^{a_\Delta} \ar@/^1pc/[drrr]^\id \ar@/_1pc/[dddr]_\id &&& \\ 
& Z \ar@<-0.5ex>[dd]\ar@<0.5ex>@{-->}[dd] \ar@<-0.5ex>@{-->}[rr]\ar@<0.5ex>[rr] & & 
X \ar@<-0.5ex>@{-->}[dd] \ar@<0.5ex>[dd]^{a_!} \\ 
& & \# &  \\
& X \ar@<-0.5ex>[rr]_{a_!'}  \ar@<0.5ex>@{-->}[rr]      & & Y      }
\xymatrix@R=1.3pc@C=1.3pc{
& &  \\ \\ \\
&\ar@{~>}[r] & \\
& &      } \qq
\xymatrix@R=1.6pc@C=1.5pc{
I \ar[dr]^\Delta \ar@/^1pc/[drrr]^\id \ar@/_1pc/[dddr]_\id &&& \\ 
& I\times I \ar[dd]\ar[rr] & & 
I \ar[dd]  \\ 
& & \pb &  \\
& I \ar[rr]     & & 1      }
\]
\epf
\begin{proposition}
If the identities over $I$ of the monoids $m$ and $m'$ commute,
then they coincide.
\end{proposition}
\pf
The left hand square below commutes by hypothesis.
\[
\xymatrix@R=1.6pc@C=1.5pc{
& X \ar@<-0.5ex>[dd]_\id\ar@<0.5ex>@{-->}[dd] \ar@<-0.5ex>@{-->}[rr]\ar@<0.5ex>[rr]^\id & & 
X \ar@<-0.5ex>@{-->}[dd] \ar@<0.5ex>[dd]^{a_\Delta} \\ 
& & \# &  \\
& X \ar@<-0.5ex>[rr]_{a_\Delta'}  \ar@<0.5ex>@{-->}[rr]      & & Z      }
\xymatrix@R=1.3pc@C=1.3pc{
& &  \\ \\
&\ar@{~>}[r] & \\ 
& &      } \q
\xymatrix@R=1.6pc@C=1.5pc{
& I \ar[dd]_\id \ar[rr]^\id & & 
I \ar[dd]^{\Delta}  \\ 
& & \pb &  \\
& I \ar[rr]_{\Delta}  & & I\times I      }
\]
\epf


\section{Connections with fibered categories}
\label{fc}

We presently show how fibered categories 
(see for instance \cite{benabou} and \cite{jacobs})
arise as a particular instance of fibered multicategories.

\begin{proposition}
\label{fc1}
Let $\I$ be a category with pullbacks.
Any fibered category $p:\M\to\I$ gives a fibered multicategory $\M = (d, p, \sp(\M,\D))$,
where $d:\D\to\I$ is the restriction of $p$ to the subcategory of cartesian arrows in $\M$,
and the special squares are the commutative squares in $\M$, with cartesian horizontal sides,
over pullbacks in $\I$.
\end{proposition}
\pf
It is easily checked that all the requirements in the definition of fibered multicategory are satisfied.
\epf

Conversely we can characterize those fibered multicategories for which $p:\M\to\I$
is a fibration.
\begin{proposition}
\label{fibchar}
Let $\M = (d, p, \sp(\M,\D))$ a fibered multicategory, and let $F:\D\to\M$ be a functor over $\I$ 
which is the identity on objects and such that for any special square 
$(f,h,a,b)$, the corresponding square $(Ff,Fh,a,b)$ in $\M$ is a pullback.
Then $p:\M\to\I$ is a fibration.
\eq
\label{fibchar2}
\xymatrix@R=1.3pc@C=1.3pc{
U \ar[dd]_b\ar@{-->}[rr]^h  & & X \ar[dd]^a \\ 
& \# & \\
V \ar@{-->}[rr]_f                 & & Y      }
\qq
\xymatrix@R=1.3pc@C=1.3pc{
& & \\ 
& \imp & \\
& &      }
\qq
\xymatrix@R=1.3pc@C=1.3pc{
U \ar[dd]_b\ar[rr]^{Fh}  & & X \ar[dd]^a \\ 
& \pb & \\
V \ar[rr]_{Ff}                 & & Y      }
\eeq

\end{proposition}
\pf
Note that the pairs $(f,Ff)$ are endomorphisms, in the sense of definition \ref{endo}.
It is enough to show that any arrow $Ff$ is $p$-cartesian.
Suppose then that $a$ is an arrow in $\M$ such that $pa = fl$ in $\I$.
We have to show that $a$ factors uniquely as $(Ff)t$ in $\M$, with $pt = l$.
Consider a reindexing $b$ of $a$ along $f$ and the pullback associated to that special square.
Since the right hand diagram below is a pullback in $\I$, we get the dotted arrow $w$.
By a $d$-lifting we get the $F w$ in the left hand diagram and the desired factorization
with $t = b(Fw)$.
\[
\xymatrix@R=1.3pc@C=1.3pc{
X \ar[dr]^{Fw} \ar@/^1pc/[drrr]^\id  &&& \\
& U \ar[dd]_b\ar[rr]^{Fh}  & & X \ar[dd]^a \\ 
&& \pb & \\
& V \ar[rr]_{Ff}                 & & Y      }
\xymatrix@R=1.3pc@C=1.3pc{
& & \\  \\  \\
& \ar@{~>}[r] & \\
& &      }
\xymatrix@R=1.3pc@C=1.3pc{
R \ar@{..>}[dr]^w \ar@/^1pc/[drrr]^\id \ar@/_1pc/[dddr]_l &&& \\
& K \ar[dd]_{pb}\ar[rr]^{h}  & & I \ar[dd]^{pa} \\ 
&& \pb & \\
& L \ar[rr]_{f}                 & & J      }
\]
As for unicity, suppose we have such a factorization $a = (Ff)t$ and consider
the map $u$ induced by the left hand pullback below:
\[
\xymatrix@R=1.3pc@C=1.3pc{
X \ar@{..>}[dr]^u \ar@/^1pc/[drrr]^\id \ar@/_1pc/[dddr]_t &&& \\
& U \ar[dd]_b\ar[rr]^{Fh}  & & X \ar[dd]^a \\ 
&& \pb & \\
& V \ar[rr]_{Ff}                 & & Y      }
\xymatrix@R=1.3pc@C=1.3pc{
& & \\  \\  \\
& \ar@{~>}[r] & \\
& &      }
\xymatrix@R=1.3pc@C=1.3pc{
R \ar[dr]^{pu} \ar@/^1pc/[drrr]^\id \ar@/_1pc/[dddr]_l &&& \\
& K \ar[dd]_{pb}\ar[rr]^{h}  & & I \ar[dd]^{pa} \\ 
&& \pb & \\
& L \ar[rr]_{f}                 & & J      }
\]
Since the right hand diagram is a pullback too, we have $pu = w$, so that $u = Fw$ 
by the unicity of the $d$-cartesian factorization $\id = (Fw)(Fh)$.
\epf

\subsection{Unary multicategories and the Grothendieck construction}
\label{gro}

As a consequence of the above proposition, we see that the free (left adjoint) way to
make a fibered multicategory into a fibered category is to freely add the arrows of $\D$
to those of $\M$ in such a way that special squares become pullback squares.
While in general the resulting fibered category has not a simple form, 
there is a well known special case. 

We say that a fibered multicategory is {\em unary} if all the maps in $\M$ are projected 
by $p:\M\to\I$ to isomorphisms.
In this case, the data of the fibered multicategory amount essentially to those of a 
pseudofunctor $\I \to \Cat$. 
In fact, if we restrict to vertical arrows, reindexing along special squares of the form
\[
\xymatrix@R=1.3pc@C=1.3pc{
U \ar@{..>}[dd]_b\ar@{-->}[rr]^f  & & X \ar[dd]^a \\ 
& \# & \\
V \ar@{-->}[rr]_f                 & & Y      }
\qq
\xymatrix@R=1.3pc@C=1.3pc{
& & \\ 
& \ar@{~>}[r]  & \\
& &      }
\qq
\xymatrix@R=1.3pc@C=1.3pc{
I \ar[dd]_\id\ar[rr]^f  & & J \ar[dd]^\id \\ 
& \pb & \\
I \ar[rr]_f                 & & J      }
\]
gives functors $f^*:\M_J \to \M_I$.
Now, in this case the well known way to freely add the arrows of $\D$ to those of $\M$
is given by the Grothendieck construction.

Thus, the equivalence between pseudofunctors and fibrations can be seen 
as the fibered version of the equivalence between unary multicategories and categories.

\subsection{Connections with sequential multicategories}
\label{seq}

If we apply the Grothendieck construction to the usual reindexing functor on $\Set_f$ 
of finite families of objects and arrows of a category $\C$,
we obtain a fibered category which, as a fibered multicategory, 
is the standard multicategory $\C\t$ of discrete cocones in $\C$
(that is, arrows in $\C\t$ are families of concurrent arrows in $\C$).

In \cite{pisani}, discrete cocone multicategories are characterized abstractly as those
multicategories with a \"central monoid", that is those in which each object carries a monoid
such that any arrow is a monoid morphism. 
There, they are named \"sequential multicategories", a somewhat unhappy choice, 
especially in the light of the present emphasis on families in place of sequences.

\begin{remark}
Note that among the consequences of the central monoid characterization
of sequential multicategories there is the well known characterization of cartesian 
monoidal categories (that is, those in which $\otimes = \times$) as those with a central comonoid. 
Indeed, a sequential multicategory is representable if and only if it is cocartesian monoidal,
which gives the dual result.
\end{remark}

As we have just mentioned, in the present framework sequential multicategories
are just the standard fibered multicategories associated to the usual family
fibration on $\Set_f$ of a category $\C$.
Thus, we are led to consider, more generally, any fibration as a sort of fibered
sequential multicategory. 
This view is also substantiated by the observation that
if $\M = (d, p, \sp(\M,\D))$ is the fibered multicategory associated to a fibration, 
then any object carries a monoid, since any section of $d:\D\to\I$ gives, via $F:\D\to\M$, 
a section of $p:\M\to\I$ (recall the definition of object and of monoid in section \ref{monin}).
Furthermore, diagram (\ref{fibchar2}) says in a sense that any map in $\M$
is \"locally" a monoid morphism. Thus, fibered categories have indeed also a sort 
of \"central monoid characterization" among fibered multicategoris.

In \cite{pisani} it is also proved that the sequential reflection of a symmetric multicategory $\M$ 
is given by $\M\otimes_{BV}1\t$, where $\otimes_{BV}$ is the Boardman-Vogt tensor product
and $1\t$ is the terminal multicategory. By the definition of the Boardman-Vogt tensor product,
$\M\otimes_{BV}1\t$ is essentially the result of adding a monoid structure to each object of $\M$, 
in such a way that any arrow of $\M$ commutes with these monoids.
This generalizes to the \"free fibration" construction on a fibered multicategory
(as sketched in section \ref{gro} above), where one adds the arrows of $\D$ to those of $\M$
in such a way that the special squares (giving reindexing) become commutative squares.


\section{Cartesian fibered multicategories}
\label{cm}

Cartesian multicategories are usually considered as a variant of Lawvere theories. 
Indeed, a (multisorted) algebraic theory can be encoded in a
cartesian multicategory $\M$, and its models (say in sets) are the cartesian functors $\M \to \Set$.
One advantage with respect to finite product categories is that in a cartesian multicategory
the products are in principle only \"virtual" , so that they have usually a simpler description.
For example, if $R$ is a ring, the Lawvere theory for $R$-modules is given by {\em matrices} with entries in $R$,
while the associated cartesian multicategory is given by {\em sequences} with entries in $R$
(see example \ref{ring}).
Indeed, the idea of a cartesian multicategory is that if it is representable 
(that is, corresponds to a monoidal category)
then it is represented by universal products (that is, the tensor product is cartesian). 
This is made precise in theorem \ref{maincart}, where it is shown that 
a fibered cartesian multicategory is stably representable if and only if it has universal products
(as defined in Section \ref{fm}). 
Indeed, both conditions are equivalent to the existence of \"algebraic products", that are defined
equationally in cartesian multicategories and then are \"absolute" with respect to cartesian functors.

There are various equivalent definitions of cartesian multicategories in the literature 
(see for instance \cite{pisani} or \cite{shulman16}).
Usually, they are considered as a sort of variation of symmetric multicategories, where
the action of bijective mappings on arrows is extended to an action of all mappings.
From our point of view, we have seen that the action of bijective mappings on the arrows 
of a symmetric multicategory is related to the {\em contravariant} reindexing of arrows 
along pullbacks (see section \ref{action}).
We will show now that the action on arrows in a cartesian multicategory is instead 
given by a sort of {\em covariant} reindexing of arrows, related to the contravariant 
one by some kind of Beck-Chevalley and Frobenius laws.

\subsection{Cartesian structures}
Given a fibered multicategory $\M = (d, p, \sp(\M,\D))$, we can consider the category $\tr(\M,\D)$,
whose objects are arrows in $\M$ and whose morphisms are triangles with the top side in $\D$,
such that their image commutes in $\I$:
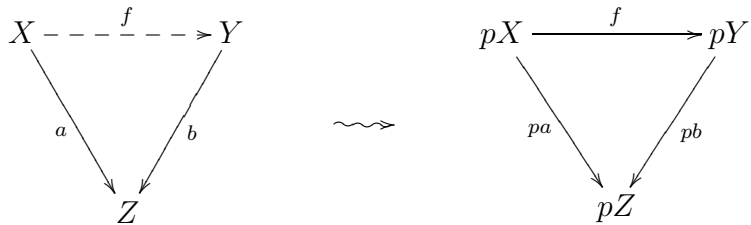
\[
\xymatrix@R=2pc@C=2pc{
X \ar[ddr]_a\ar@{-->}[rr]^f  & & Y \ar[ddl]^b \\ 
& &  & \ar@{~>}[r] & \\
             & Z &      }    \qq
\xymatrix@R=1.8pc@C=1.8pc{
pX \ar[ddr]_{pa}\ar[rr]^f  & & pY \ar[ddl]^{pb} \\ 
&& \\
             & pZ &      }
\]
Composition and identities are the obvious ones.

\subsection{Construction C}
Let $\tr(d)$ be the comma category $(d,\delta_\I)$, where $\delta_\I:\obj\I\to\I$ is the discrete category inclusion.
Thus an object in $\tr(d)$ is a pair of objects $I\in \I$ and $X\in\D$ along with an arrow $l:d X \to I$, 
while an arrow from $l:d X \to I$ to $h:d Y \to I'$ exists only if $I=I'$ and is an arrow $a:X\to Y$ in $\D$
such that $h(da) = l$.
There is an obvious functor $\tr(\M,\D) \to \tr(d)$.
\begin{definition} \label{defcart}
A {\em cartesian structure} on a fibered multicategory $\M$ consists of a subcategory of \"special triangles"
$\sptr(\M,\D) \subseteq \tr(\M,\D)$, such that:
\begin{enumerate}
\item
Composing a special triangle $(a,b,f)$ with an arrow $c:Z\to W$ in $\M$ 
(where $Z$ is the common codomain of $a$ and $b$)
one gets a special triangle $(ca,cb,f)$.
\item
The obvious functor $\sptr(\M,\D) \to \tr(d)$ is a discrete opfibration.
\item
Pasting a special triangle with a special square over it, one gets a special triangle.
\item
Pulling back a special triangle along special squares, one gets a special triangle.
\end{enumerate}
A {\em cartesian morphism} $\M\to\M'$ of cartesian multicategories, is a morphism 
of multicategories that sends special triangles in $\M$ to special triangles in $\M'$.
\end{definition}

Thus, special triangles and special squares are related by what we will refer respectively 
as the Frobenius and Beck-Chevalley laws. 
The first one says that the following pasting is still special:
\eq
\label{FR}
\xymatrix@R=2pc@C=2pc{
X' \ar[dd]_a\ar@{-->}[rr]^{f'}  & & Y' \ar[dd]^b \\ 
& \# &  \\
X \ar@/_/[ddr]_c\ar@{-->}[rr]^f  & & Y \ar@/^/[ddl]^d \\ 
& \# &  \\
& Z &      }   
\eeq
The second one says that, in the prism below, if the top face commutes in $\D$,
 the front and back faces are special squares and the right face is a special triangle, 
 then the left face is also a special triangle:
\eq
\label{BC}
\xymatrix@R=1.5pc@C=1.5pc{
&& Y' \ar[dddl] \ar@{-->}[rrrr] &&&& Y\ar[dddl] \\
X' \ar[ddr]\ar@{-->}[rru]^{f'}  \ar@{-->}[rrrr] &&&& X\ar[ddr]_a\ar@{-->}[rru]^f  \\ 
&&  \\
& Z' \ar@{-->}[rrrr]_g &&&& Z &      }   
\eeq
Note that the top square is projected to a pullback in $\I$ (and so is itself a pullback in $\D$), 
since the front and the back squares, being special, are projected to pullbacks.

Condition (2) says roughly that maps in $\M$ can be covariantly reindexed along commutative triangles:
given an arrow $a:X\to Z$ in $\M$, a commutative triangle $hf = pa$ in $\I$ and a $d$-lifting of $f$,
there is a unique extension (by a map $b$) of the data to a special triangle as illustrated below:
\eq \label{cov}
\xymatrix@R=2pc@C=2pc{
X \ar[ddr]_a\ar@{-->}[rr]^f  & & Y \ar@{..>}[ddl]^b \\ 
& \# &  & \ar@{~>}[r] & \\
             & Z &      }    \qq
\xymatrix@R=1.8pc@C=1.8pc{
pX \ar[ddr]_{pa}\ar[rr]^f  & & I  \ar[ddl]^h \\ 
&& \\
             & pZ &      }
\eeq

We sometimes will denote such a unique map $b$ by $f_!a$. 
Note however that, as for contravariant reindexing, this is somewhat imprecise, 
since $b$ depends also on the chosen lifting of $f$.
Anyway, by doing so we find that the Frobenius law has the form 
\[ f'_! (c(f^*b)) = (f_!c)b \]
which is indeed very similar to the Frobenius reciprocity law.

Similarly, the Beck-Chevalley law becomes
\[ g^*(f_!a) = f'_!(g^*a) \]
and says roughly that (contravariant) reindexing preserves special triangles.

\begin{remark}
\label{defcart2}
One can define fibered cartesian multicategories by using \"covariant special squares" 
instead of special triangles, as we sketch below leaving the details to the reader.
Covariant special squares have the top side in $\D$ and lay over commutative squares in $\I$.
Furthermore they form a discrete opfibration over them in the obvious sense, 
as depicted in the diagram below:
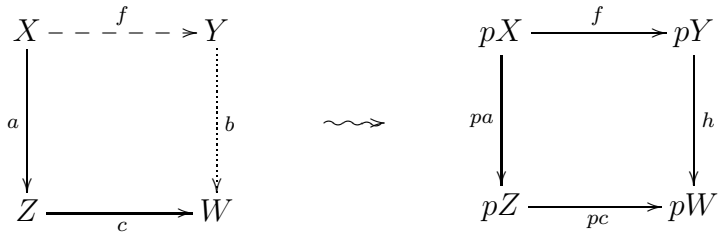
\[
\xymatrix@R=2pc@C=2pc{
X \ar[dd]_a\ar@{-->}[rr]^f  & & Y \ar@{..>}[dd]^b \\ 
& &  & \ar@{~>}[r] & \\
Z  \ar[rr]_c         & &     W       }    \qq
\xymatrix@R=1.8pc@C=1.8pc{
pX \ar[dd]_{pa}\ar[rr]^f  & & pY \ar[dd]^h \\ 
&& \\
pZ \ar[rr]_{pc} & & pW    }
\]
Then condition (1) of definition \ref{defcart} becomes redundant, since it is a particular instance
of horizontal pasting (that is, composition) of covariant special squares, 
namely when the top side of the second square is the identity.
The Frobenius condition is a sort of vertical pasting of a contravariant and a covariant special square.
The Beck-Chevalley condition is now depicted by the cube below, in place of (\ref{BC}):
if the top face commutes in $\D$, the front, back and bottom faces are special squares 
and the right face is a covariant special square, then the left face is also a covariant special square.
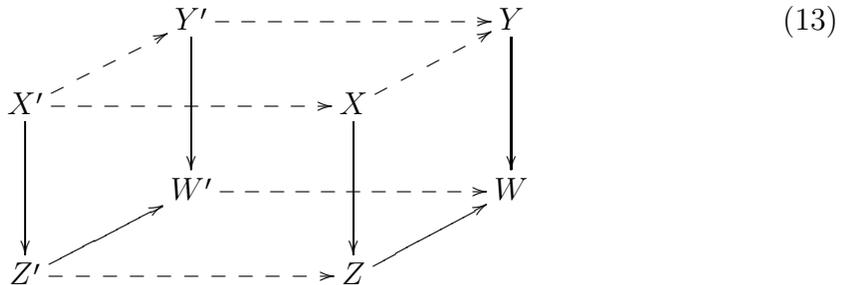
\eq
\label{BC2}
\xymatrix@R=1.5pc@C=1.5pc{
&&& Y' \ar[dd] \ar@{-->}[rrrr] &&&& Y\ar[dd] \\
&X' \ar[dd]\ar@{-->}[rru]  \ar@{-->}[rrrr] &&&& X\ar[dd]\ar@{-->}[rru]  \\ 
&& & W' \ar@{-->}[rrrr] &&&& W  \\
& Z' \ar@{-->}[rrrr] \ar[rru] &&&& Z \ar[rru] &      }   
\eeq
It is then easy to see that the two definitions are equivalent. 
\end{remark}

\subsection{Standard cartesian multicategories}

Le us show how in the standard case (see section \ref{standard}) we find again 
(the \"non-skeletal" version of) classical cartesian multicategories.
First, note that if $pZ$ is terminal in $\Set_f$, the right hand triangle in (\ref{cov}) 
always commutes and condition (2) of definition \ref{defcart} becomes:
for any $d$-arrow $X\to Y$ over $f:I\to J$ (that is, $Y = y_j \,\, (j\in J)$ is a family of objects in $\M_0$, 
and $X = x_i = y_{f(i)} \,\, (i\in I)$ is the family obtained by reindexing $Y$ according to $f$) 
and for any arrow $a:X\to Z$ in $\M$, there is a uniquely determined covariant reindexing $f_! a : Y\to Z$.
If $pZ$ is not terminal, since the right hand triangle in (\ref{cov}) commutes, 
$f$ is in fact a \"sum" of mappings and, by the Beck-Chevalley condition, the covariant reindexing along $f$
reduces to the covariant reindexing of each single arrow in the family.
Furthermore, condition (1) of definition \ref{defcart} says that reindexing a family of arrows and then 
composing it with an arrow, gives the same result as first composing and then reindexing the resulting arrow. 

Since $\sptr(\M,\D)$ is a subcategory of $\tr(\M,\D)$, the covariant reindexing is functorial
and, lastly, the Frobenius law gives, in a synthetical form, the (otherwise somewhat cumbersome) 
\"block reindexing" condition, concerning the composition of a family of arrows with a reindexed arrow. 

\begin{example}
If $\M$ is associated to a finite product category $\C$, the idea is that one can duplicate or delete the variables.
For instance, given a map $t : A\tm B \tm A \to D$ in $\Set$, one gets another
map $f_! t: B\tm C \tm A \to D$ by
\[ 
f_! t(b,c,a) = t(a,b,a) \,\,.
\]
To be precise, if $Y$ is the family of sets $\{1,2,3\} \to\obj\Set$ given by $1\mapsto B,2\mapsto C,3\mapsto A$ 
and $f:\{1,2,3\} \to\{1,2,3\}$ is $1\mapsto 3,2\mapsto 1,3\mapsto 3$, then the domain $X=f^*Y$ of $t$
is the family of sets $\{1,2,3\}\to\obj\Set$ given by $1\mapsto A,2\mapsto B,3\mapsto A$;
$f_!t$ is obtained by precomposing $t$ with the obvious map $\ov f : B\tm C \tm A \to A\tm B\tm A$.
As particular instances of $\ov f$ we find projections and diagonals.
\end{example}

\begin{example}
\label{ring}
If $\C$ is any additive category (or also any category enriched in commutative monoids), 
then there is a natural cartesian structure on the corresponding sequential multicategory 
$\C\t$ (see \cite{pisani}). 
For instance, a ring $R$ gives a cartesian multicategory $\M_R$ with just one object 
and whose arrows are families of element of $R$.
If $t:\{1,2,3\} \to R$ is the arrow $1\mapsto b,2\mapsto c,3\mapsto a$ in $\M_R$, 
and $f:\{1,2,3\} \to\{1,2,3\}$ is $1\mapsto 3,2\mapsto 1,3\mapsto 3$ as in the example above,
then $f_!t$ is obtained by summing or inserting $0$'s: $1\mapsto c,2\mapsto 0,3\mapsto b+a$.
\end{example}

\subsection{Coherence of the two reindexing}

Now, one may wonder whether the two indexing are coherent, when both are possible.
This is in fact guaranteed by the Beck-Chevalley law, as shown in the following proposition.
Given a cartesian fibered multicategory $\M$, a commutative triangle in $\I$ whose top side is 
an {\em isomorphism}, and a lifting $a:X\to Y$ of its left side in $\M$,
we can get not only the covariant reindexing $f_!a$, but also the contravariant
reindexing $g^*a$ along the inverse map $g = f^{-1}$:
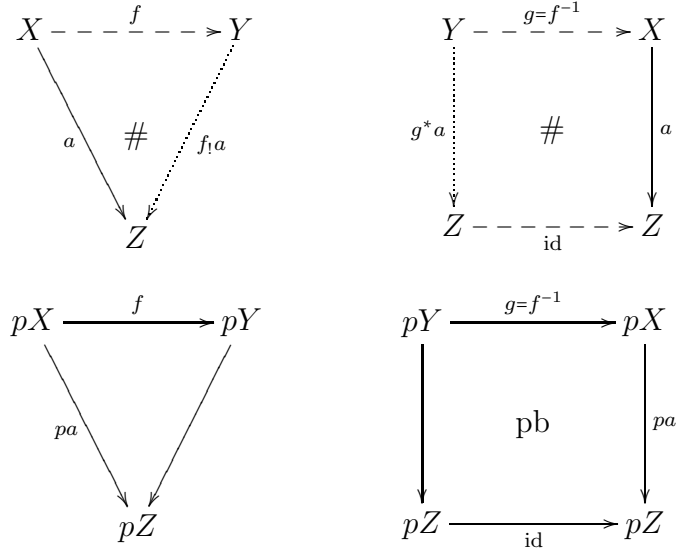
\[
\xymatrix@R=2pc@C=2pc{
X \ar[ddr]_a\ar@{-->}[rr]^f  & & Y \ar@{..>}[ddl]^{f_!a} \\ 
& \# & & \\
& Z &      }    \qq
\xymatrix@R=1.8pc@C=1.8pc{
Y \ar@{..>}[dd]_{g^*a}\ar@{-->}[rr]^{g=f^{-1}}  & & X \ar[dd]^a \\ 
& \# & \\
Z \ar@{-->}[rr]_\id                 & & Z      }             
\]
\[
\xymatrix@R=2.3pc@C=1.5pc{
pX \ar[ddr]_{pa}\ar[rr]^f  & & pY  \ar[ddl] \\ 
&& \\
             & pZ &      } \qq\qq
\xymatrix@R=1.8pc@C=1.8pc{
pY \ar[dd]\ar[rr]^{g=f^{-1}}  & & pX \ar[dd]^{pa} \\
& {\rm pb} & \\
pZ  \ar[rr]_\id                 & & pZ       }             
\]

\begin{prop}
In the above situation, the two reindexing coincide: $f_!a = g^*a$.
\end{prop}
\pf
Consider the following prism, where $g = f^{-1}$ (so that the top square commutes)
and where we assume that the right triangle and the front and back squares are special.
Then, by the Beck-Chevalley condition, the left triangle is also special, so that $t = g^*a$.
But since the back square is special, also $t = f_!a$.
\[
\xymatrix@R=1.5pc@C=1.5pc{
&& Y \ar[dddl]^(.6)t \ar@{-->}[rrrr]^\id &&&& Y\ar[dddl]^{f_!a} \\
Y \ar[ddr]_{g^*a}\ar@{-->}[rru]^\id  \ar@{-->}[rrrr]^(.6)g &&&& X\ar[ddr]_a\ar@{-->}[rru]^(.5)f  \\ 
&&  \\
& Z \ar@{-->}[rrrr]_\id &&&& Z &      }   
\]
\epf

\subsection{Algebraic products}

We have defined universal products in section \ref{unipro}.
Now we show how in a fibered cartesian multicategory there is an \"algebraic" notion
of product, which subsumes both universal products and opcartesian arrows.

If $X\in\M$ and $f:pX\to J$ is a map in $\I$, an {\em algebraic product} for $X$ along $f$ 
is an object $P\in\M$ over $J$ along with a vertical map $\pi:f^*P\to X$ and a map $u:X\to P$ over $f$
\eq \label{algpr1}
\xymatrix@R=1.3pc@C=1.3pc{
f^*P \ar[dd]_\pi\ar@{-->}[ddrr]^f  & &  \\ 
& & \\
X    \ar[rr]_u            & & P      } \q
\xymatrix@R=1.3pc@C=1.3pc{
& &  \\ 
&\ar@{~>}[r] & \\
& &      } \qq
\xymatrix@R=1.3pc@C=1.6pc{
pX \ar[dd]_\id\ar[ddrr]^f  & &  \\ 
& &  \\
pX  \ar[rr]_f               & & J      } 
\eeq
such that the following are both special triangles:
\eq. \label{algpr2}
\xymatrix@R=2pc@C=0.5pc{
f^*P \ar[d]_\pi\ar@{-->}[rr]^f  & & P \ar@/^1pc/[ddl]^\id \\ 
X \ar@/_/[dr]_u & \# & & \\
& P &      }    \qq \qq \qq
\xymatrix@R=2pc@C=0.5pc{
X \ar@/_1pc/[ddr]_\id\ar@{-->}[rr]^\Delta  & & h^*X \ar[d]^{f^*u} \\ 
& \# & f^*P \ar@/^/[dl]^\pi & \\
& X &      }    
\eeq
where the $d$-map $\Delta$ in the right hand triangle is a $d$-lifting of 
the diagonal of the pullback of $f$ along itself in $\I$,
such that $h\Delta = \id$ also in $\D$:
\[
\xymatrix@R=1.2pc@C=1pc{
X  \ar@{-->}[dr]^\Delta \ar@{-->}@/^0.9pc/[drrr]^\id &&& \\
& h^*X \ar[dd]_{f^*u}\ar@{-->}[rr]^h  & & X \ar[dd]^u \\
& & \# & \\
& f^*P  \ar@{-->}[rr]_f                 & & P       }     \qq
\xymatrix@R=1.3pc@C=1.3pc{
& &  \\ \\
&\ar@{~>}[r] & \\
& &      } \qq
\xymatrix@R=1pc@C=1pc{
pX  \ar[dr]^\Delta  \ar@/_0.7pc/[dddr]_\id \ar@/^0.7pc/[drrr]^\id &&& \\
& K \ar[dd]^l\ar[rr]^h  & & pX \ar[dd]^f \\
& & \pb & \\
& pX  \ar[rr]_f                 & & J       }             
\]

Now we come to the main result of this section.

\begin{theorem}
\label{maincart}
For a fibered cartesian multicategory $\M$, the following are equivalent:
\begin{enumerate}
\item
$\M$ has algebraic products. {\q \em (AP)}
\item
$\M$ has universal products. {\q \em (UP)}
\item
$\M$ is stably representable. {\q \em (SR)}
\end{enumerate}
\end{theorem}
\pf 
Let us show that (AP) implies (UP). Suppose then that $(\pi,u)$ is an algebraic product
for $X$ along $f$ as in (\ref{algpr1}). 
We want to show that $\pi$ is a universal product for $X$ along $f$ as in (\ref{unipr}):
given a pullback $hf' = fh'$ in $\I$ with a $d$-lifting or its top side $f'$,
any arrow $\rho: f'^*Q\to X$ with $p\rho = h'$, 
should factor uniquely as $\rho = \pi(f^*t)$ for a unique special square: 
\[
\xymatrix@R=1.3pc@C=1.3pc{
& f'^*Q \ar@{..>}[ddl]_{f^*t} \ar@/^/[ddddl]^(0.3)\rho \ar@{-->}[ddrr]^{f'}&& \\ \\
f^*P \ar[dd]_\pi\ar@{-->}[ddrr]^f  & & \# & Q \ar@{..>}[ddl]^t \\ 
& &  \\
X                 & & P      } 
\xymatrix@R=2pc@C=2pc{
 \\ \\ & \ar@{~>}[r] &      } \qq
\xymatrix@R=1.3pc@C=1.3pc{
& L \ar[ddl]_{h'} \ar@/^/[ddddl]^(0.3){h'} \ar[ddrr]^{f'}&& \\ \\
pX \ar[dd]_\id\ar[ddrr]^f  & & \pb & K \ar[ddl]^h \\ 
& &  \\
pX                 & & J      } 
\]

Supposing that such a factorization exists, by the Frobenius law applied to
the first of (\ref{algpr2}), we see that $t = f'_!(u\rho)$:
\[
\xymatrix@R=2pc@C=1.5pc{
f'^*Q \ar[dd]_{f^*t}\ar@{-->}[rr]^{f'}  & & Q \ar[dd]^t \\ 
& \# & \\
f^*P \ar[d]_\pi\ar@{-->}[rr]^f  & & P \ar@/^0.8pc/[ddl]^\id \\ 
X \ar@/_/[dr]_u & \# & & \\
& P &      }    \qq \qq \qq
\xymatrix@R=2pc@C=1.5pc{
f'^*Q \ar[dd]_\rho\ar@{-->}[rr]^{f'}  & & Q \ar@/^1.3pc/[ddddl]^t \\ \\
X \ar@/_/[ddr]_u & \# & & \\ \\
& P &      } 
 \]
This proves unicity. To show that such a $t$ gives indeed the desired factorization,
consider the diagram below:
\[
\xymatrix@R=1pc@C=1pc{
&&& f'^*Q\ar@{-->}[rr]^{f'}\ar@/^1pc/[dddddl]|(.4){f^*t} & & Q  \ar@/^1pc/[dddddl]^t  \\
f'^*Q\ar@{-->}@/^0.8pc/[urrr]^\id\ar@{-->}[rr]^\Delta \ar[dd]_\rho && 
W \ar@{-->}[rr]\ar@{-->}[ur]\ar[dd]_{h^*\rho} && f'^*Q \ar[dd]_\rho\ar@{-->}[ur]  \\ 
&&&& & & & \\
X\ar@{-->}[rr]^\Delta\ar@/_1pc/[rrdddd]_\id && h^*X \ar@{-->}[rr]^h\ar[dd]_{f^*u} && X \ar[dd]_u & & & \\ \\
&& f^*P \ar@{-->}[rr]^f\ar[dd]_\pi && P &&&    \\ \\  
&& X }
\]
The right-hand triangle (with $t$ as a side) is special by hypothesis, 
so by Beck-Chevalley also the triangle with $f^*t$ as a side is special.
The lower left-hand triangle is special by the second of (\ref{algpr2}) so that,
by Frobenius, it is special also its pasting with the left-hand special square.
Now, by composition, we get a special triangle with an identity top side,
so that $\rho = \pi(f^*t)$, as desired.

Next, we prove that (AP) implies (SR). First let us show that
in an algebraic product $(\pi,u)$ along $f:pX\to J$, $u$ is opcartesian in $\M$ over $f$.
So let $v:X\to Q$ and let $g: J\to pQ$ be such that $gf = pv$;
then there should be a unique map $t:P\to Q$ over $g$, such that $tu = v$:
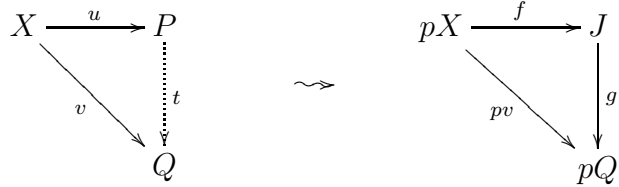
\[
\xymatrix@R=1.3pc@C=1.3pc{
X  \ar[ddrr]_v  \ar[rr]^u    & & P \ar@{..>}[dd]^t  \\
& & \\
& & Q  } \q
\xymatrix@R=1.3pc@C=1.3pc{
& &  \\ 
&\ar@{~>}[r] & \\
& &      } \qq
\xymatrix@R=1.3pc@C=1.3pc{
pX  \ar[ddrr]_{pv}  \ar[rr]^f    & & J \ar[dd]^g  \\
& & \\
& & pQ  }
\]
Supposing that such a factorization exists, by the condition (1) in the definition
of cartesian multicategory applied to the first of (\ref{algpr2}), we see that $t = f_!(v\pi)$:
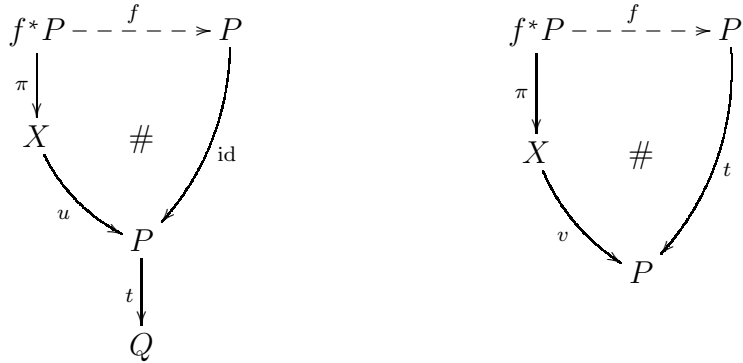
\[
\xymatrix@R=2pc@C=1.5pc{
f^*P \ar[d]_\pi\ar@{-->}[rr]^f  & & P \ar@/^0.8pc/[ddl]^\id \\ 
X \ar@/_/[dr]_u & \# & & \\
& P \ar[d]_t &   \\  
& Q &   }    \qq \qq \qq
\xymatrix@R=1pc@C=1.5pc{
f^*P \ar[dd]_\pi\ar@{-->}[rr]^f  & & P \ar@/^1pc/[ddddl]^t \\ \\
X \ar@/_/[ddr]_v & \# & & \\ \\
& P &      } 
 \]
This proves unicity. To show that such a $t$ gives indeed the desired factorization,
consider the diagram below:
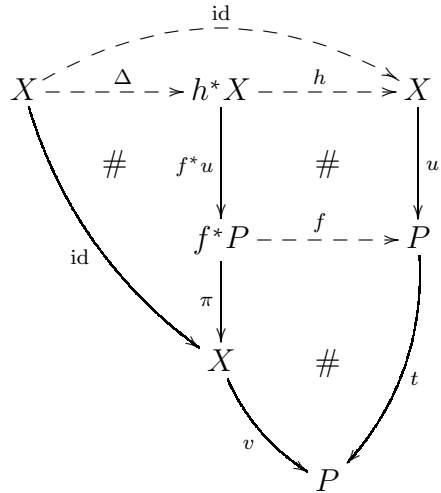
\[
\xymatrix@R=1pc@C=1.5pc{
X \ar@{-->}[rr]^\Delta \ar@/_1pc/[ddddrr]_\id \ar@/^2pc/@{-->}[rrrr]^\id 
&& h^*X \ar[dd]_{f^*u}\ar@{-->}[rr]^h  & & X \ar[dd]^u \\ 
& \# & &\# & \\
&& f^*P \ar[dd]_\pi\ar@{-->}[rr]^f  & & P \ar@/^1pc/[ddddl]^t \\ \\
&& X \ar@/_/[ddr]_v & \# & & \\ \\
&& & P &      } 
 \]
The outer triangle is special (by pasting and compositions) and has an identity top side,
so that $v = tu$, as desired.
To show that $u$ is {\em stably} opcartesian, note that the notion of algebraic product is 
stable with respect to reindexing: 
if $(\pi,u)$ is an algebraic product along $f:pX\to J$ and $l:L\to J$ is a map in $\I$,
then reindexing the (\ref{algpr1}) along $l$ one gets an algebraic product $(\pi',u')$
along $f'$ (the pullback of $f$ along $l$).
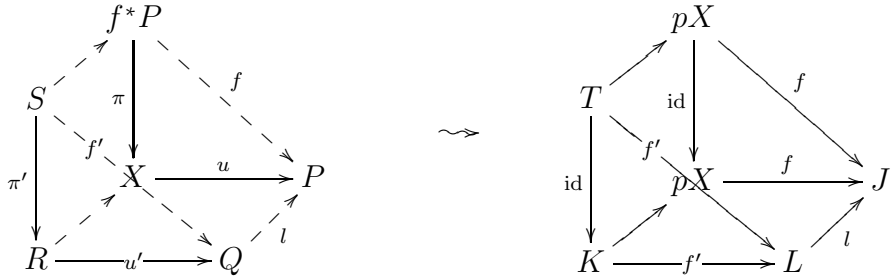
\[
\xymatrix@R=1.3pc@C=1.3pc{
& f^*P \ar[dd]_\pi\ar@{-->}[ddrr]^f  & &  \\ 
S \ar@{-->}[ur]\ar[dd]_{\pi'} \ar@{-->}[ddrr]|(.3){f'} & & \\
& X    \ar[rr]^u            & & P   \\
R \ar@{-->}[ur]\ar[rr]|{u'} && Q \ar@{-->}[ur]_l    } \q
\xymatrix@R=1.3pc@C=1.3pc{
& &  \\  \\
&\ar@{~>}[r] & \\
& &      } \qq
\xymatrix@R=1.3pc@C=1.6pc{
& pX \ar[dd]_\id\ar[ddrr]^f  & &  \\ 
T \ar[ur]\ar[dd]_\id\ar[ddrr]|(.3){f'} & &  \\
& pX  \ar[rr]^f               & & J  \\
K \ar[ur]\ar[rr]|{f'} && L \ar[ur]_l   } 
\]
Indeed, since (contravariant) reindexing preserves identity, composition, special squares
and special triangles, it preserves the \"equations" (\ref{algpr2}) as well. 

Next, we prove that (UP) implies (AP). Suppose then that $\pi:f^*P\to X$
is a universal product for $X$ along $f$:
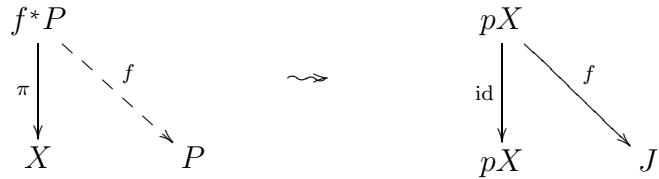
\[
\xymatrix@R=1.3pc@C=1.3pc{
f^*P \ar[dd]_\pi\ar@{-->}[ddrr]^f  & &  \\ 
& &  \\
X                 & & P      } 
\xymatrix@R=1.3pc@C=1.3pc{
& &  \\ 
&\ar@{~>}[r] & \\
& &      } \qq\qq
\xymatrix@R=1.3pc@C=1.3pc{
pX \ar[dd]_\id\ar[ddrr]^f  & &  \\ 
& &  \\
pX                 & & J      } 
\]
We show that $\pi$ is part of an algebraic product $(\pi,u)$. 
The map $u:X\to P$ is obtained by exploiting the universal property of $\pi$
(in fact, the \"existence" part).
It is the (unique) map such that $\pi(f^*u) = \Delta_!\id$, that is such that
both the square and the triangle in the diagram below are special:
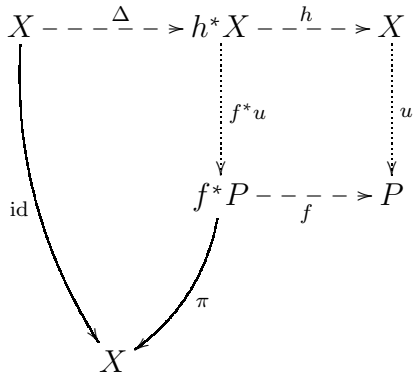
\[
\xymatrix@R=4pc@C=1.5pc{
X \ar@/_1pc/[ddr]_\id \ar@{-->}[rr]^\Delta  & & h^*X \ar@{..>}[d]^{f^*u}\ar@{-->}[rr]^h && X \ar@{..>}[d]^u \\ 
&  & f^*P \ar@/^.8pc/[dl]^\pi \ar@{-->}[rr]_f & & P \\
& X &      }    
\]
To show that also the first of equations (\ref{algpr2}) holds true, we exploit again
the universality of $\pi$ (but now the \"uniqueness" part).
We want to show that in the right hand special triangle of the diagram below, $i$ is the identity:
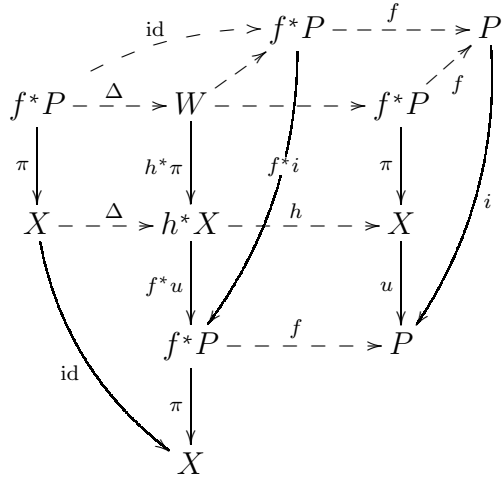
\[
\xymatrix@R=1pc@C=1pc{
&&& f^*P\ar@{-->}[rr]^{f}\ar@/^1pc/[dddddl]|(.4){f^*i} & & P  \ar@/^1pc/[dddddl]^i  \\
f^*P\ar@{-->}@/^0.8pc/[urrr]^\id\ar@{-->}[rr]^\Delta \ar[dd]_\pi && 
W \ar@{-->}[rr]\ar@{-->}[ur]\ar[dd]_{h^*\pi} && f^*P \ar[dd]_\pi\ar@{-->}[ur]_f  \\ 
&&&& & & & \\
X\ar@{-->}[rr]^\Delta\ar@/_1pc/[rrdddd]_\id && h^*X \ar@{-->}[rr]^h\ar[dd]_{f^*u} && X \ar[dd]_u & & & \\ \\
&& f^*P \ar@{-->}[rr]^f\ar[dd]_\pi && P &&&    \\ \\  
&& X }
\]
By the Beck-Chevalley and Frobenius laws, the left hand triangle with left side $\id_X\pi$ 
and right side $\pi(f^*i)$ is also special.
Since its top side is the identity, we get $\pi(f^*i) = \pi$. So $i$ and $\id_P$
both give a factorization of $\pi$ through $\pi$, and we get $i = \id_P$ as desired.

Lastly, we prove that (SR) implies (AP). 
Suppose then that $u:X\to P$ is a stably opcartesian arrow over $f$.
We show that $u$ is part of an algebraic product $(\pi,u)$. 
The map $\pi:f^*P\to X$ is obtained by exploiting the universal property of $u$
(in fact, the \"existence" part).
It is the (unique) map such that $\pi(f^*u) = \Delta_!\id$, that is 
such that the triangle in the diagram below is special:
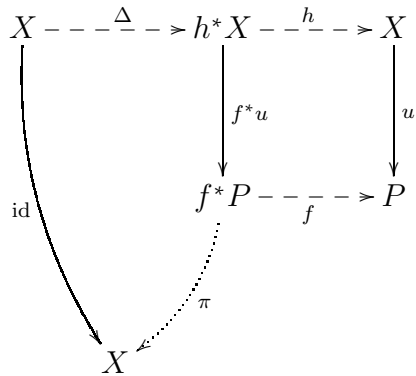
\[
\xymatrix@R=4pc@C=1.5pc{
X \ar@/_1pc/[ddr]_\id \ar@{-->}[rr]^\Delta  & & h^*X \ar[d]^{f^*u}\ar@{-->}[rr]^h && X \ar[d]^u \\ 
&  & f^*P \ar@/^.8pc/@{..>}[dl]^\pi \ar@{-->}[rr]_f & & P \\
& X &      }    
\]
To show that also the first of equations (\ref{algpr2}) holds true, we exploit again
the universality of $u$ (but now the \"uniqueness" part).
We want to show that in the bottom special triangle of the diagram below, $i$ is the identity:
\[
\xymatrix@R=1pc@C=1.5pc{
X \ar@{-->}[rr]^\Delta \ar@/_1pc/[ddddrr]_\id \ar@/^2pc/@{-->}[rrrr]^\id 
&& h^*X \ar[dd]_{f^*u}\ar@{-->}[rr]^h  & & X \ar[dd]^u \\ 
& \# & &\# & \\
&& f^*P \ar[dd]_\pi\ar@{-->}[rr]^f  & & P \ar@/^1pc/[ddddl]^i \\ \\
&& X \ar@/_/[ddr]_u & \# & & \\ \\
&& & P &      } 
 \]
The outer triangle is special (by pasting and compositions) and has an identity top side,
so that $iu = u = \id_Pu$. Thus, $i = \id_P$ as desired.
\epf

\begin{remark}
The present proof of theorem \ref{maincart} follows, in the more abstract setting of fibered
multicategories, the proof given in \cite{pisani} for the classical case.
A similar proof is sketched, in the context of LNL polycategories, in \cite{shulman21}, 
while in \cite{garner} it is shown how the result is a consequence of an elegant general 
theorem on absolute enriched colimits. 

\end{remark}

\subsection{Cartesian fibered categories}

If $\M$ is the fibered multicategory associated to a fibered category $p:\M\to\I$ 
(see proposition \ref{fc1}) then it is stably representable if and only if the fibration $p$ has sums 
(in the sense of fibered categories, see for instance \cite{benabou}).
Indeed, since in this case the special squares are commutative squares in $\M$ 
with cartesian horizontal arrows over pullbacks in $\I$,
the stability condition of definition \ref{rep} becomes the Beck-Chevalley condition for bifibrations.
Similarly, $\M$ has universal products if and only if the fibration $p$ has products, 
that is the opposite fibration has sums.

Thus, if we say that a fibered category is cartesian if the associated fibered multicategory is so,
as a corollary of theorem \ref{maincart} we get:
\begin{proposition}
A cartesian fibered category $\M$ has sums if and only if it has products and, 
if this is the case, they coincide and are absolute (with respect to the cartesian morphisms
of definition \ref{defcart}).
\end{proposition}
\epf

Recall that the usual family fibration of a category $\C$ on $\Set_f$ corresponds,
as a fibered multicategory, to the standard sequential multicategory $\C\t$
(see section \ref{seq}).
In this case, as shown in \cite{pisani}, a cartesian structure on $\C\t$ amounts to an 
enrichment of $\C$ in commutative monoids (see also example \ref{ring}).
Thus we find again the following fact, which is well-known in the case of abelian categories
\begin{proposition}
A category $\C$, enriched in commutative monoids, has finite sums if and only if 
it has finite products and, if this is the case, they coincide and are absolute.
\end{proposition}
\epf

To conclude, we present a partial generalization of the above mentioned correspondence 
between cartesian structures on $\C\t$ and enrichments of $\C$ in commutative monoids.
First, recall from section \ref{monin} that an object in a fibered multicategory is a section of $d:\D\to\I$;
in the case of fibered categories it is then a section $s:\I\to\M$ of $p:\M\to\I$ such that $sf:sI\to sJ$ 
is cartesian for any $f:I\to J$ in $\I$. Given two objects $s$ and $t$ in $\M$, let the \"fibered hom-set"
$p(s,t):\M(s,t)\to\I$ be the mapping such that for any $I\in\I$, $\M(s,t)_I$ is the set of vertical arrows $sI\to tI$.

\begin{theorem}
The fibered hom-set $p(s,t):\M(s,t)\to \I$ of any pair of objects in a cartesian fibered category over $\I$ 
has naturally the structure of a fibered monoid over $\I$ (see definition \ref{fibmon}).
\end{theorem}
\pf
The contravariant reindexing extends $p(s,t):\M(s,t)\to \I$ to a discrete fibration, while the covariant
reindexing extends it to a discrete opfibration.
\[
\xymatrix@R=1.5pc@C=1.5pc{
sI \ar@{..>}[dd]_{f^*a}\ar@{-->}[rr]^{sf}  & & sJ \ar[dd]^a \\ 
& \# & & \ar@{~>}[r] & \\
tI  \ar@{-->}[rr]_{tf}         & &     tJ       }    \qq
\xymatrix@R=2pc@C=2pc{
I \ar[dd]_\id\ar[rr]^f  & & J \ar[dd]^\id \\ 
&& \\
I \ar[rr]_f & & J    }
\]
\[
\xymatrix@R=1.5pc@C=1.5pc{
sI \ar[dd]_a\ar@{-->}[rr]^{sf}  & & sJ \ar@{..>}[dd]^{f_!b} \\ 
& \# & & \ar@{~>}[r] & \\
tI  \ar@{-->}[rr]_{tf}          & &     tJ       }    \qq
\xymatrix@R=2pc@C=2pc{
I \ar[dd]_\id\ar[rr]^f  & & J \ar[dd]^\id \\ 
&& \\
I \ar[rr]_f & & J    }
\]
Now, the condition for $p(s,t):\M(s,t)\to \I$ to be a fibered monoid becomes exactly
the alternative form (\ref{BC2}) of the Beck-Chevalley condition.
\epf

%



\begin{refs}


\bibitem[Beilinson \& Drinfeld, 2004]{beilinson} Beilinson and Drinfeld (2004), {\em Chiral algebras}, 
American Mathematical Society Colloquium Publications, 51.

\bibitem[Cruttwell \& Shulman, 2010]{shulman} G. Cruttwell and M. Shulman (2010), 
A unified framework for generalized multicategories, {\em Theory and Appl. of Cat.} {\bf 24}, 580-655.

\bibitem[Garner, 2014]{garner} R. Garner (2014), Diagrammatic characterisation of enriched absolute colimits, 
{\em Theory and Appl. of Cat.} {\bf 29}, 775–780.



\bibitem[Hermida, 2000]{hermida} C. Hermida (2000), Representable multicategories, 
{\em Advances in Math.}, {\bf 151}, 164-225.

\bibitem[Hermida, 2004]{hermida2} C. Hermida (2004), {\em Fibrations for abstract multicategories}, 
in: Galois theory, Hopf algebras, and semiabelian categories, Fields inst. comm. AMS, 281-293.

\bibitem[Jacobs, 1999]{jacobs} B. Jacobs (1999), {\em Categorical Logic and Type Theory}, Studies in Logic and the Foundations of Mathematics 141, North Holland, Elsevier. 



\bibitem[Leinster, 2003]{leinster} T. Leinster (2003), {\em Higher operads, higher categories}, Cambridge University Press, math.CT/0305049.

\bibitem[May \& Thomason, 1978]{may} J. P. May and R. Thomason (1978), The uniqueness of infinite loop space machines,
 {\em Topology} {\bf 17}, 215-224.
 
\bibitem[Pisani, 2014]{pisani} C. Pisani (2014), Sequential multicategories, {\em Theory and Appl. of Cat.} {\bf 29}, 496-541.

\bibitem[Shulman, 2016]{shulman16} M. Shulman (2016), {\em Categorical logic from a categorical point of view}, 
draft for AARMS Summer School, available on GitHub.

\bibitem[Shulman, 2021]{shulman21} M. Shulman (2021), LNL polycategories and doctrines of linear logic,
preprint available on arXiv.org.

\bibitem[Streicher, 2018]{benabou} T. Streicher (2018), Fibered Categories a la Jean Benabou,
preprint available on arXiv.org.

\end{refs}

\end{document}